\documentclass[12pt,english,a4paper]{smfart}

\usepackage[T1]{fontenc}
\usepackage{lmodern}
\usepackage{smfthm}
\usepackage[headings]{fullpage}
\usepackage{amssymb,amscd}

\newcommand{\Qp}{\mathbf{Q}_p}
\newcommand{\Cp}{\mathbf{C}_p}

\newcommand{\ZZ}{\mathbf{Z}}
\newcommand{\FF}{\mathbf{F}}
\newcommand{\QQ}{\mathbf{Q}}
\newcommand{\OO}{\mathcal{O}}
\newcommand{\MM}{\mathfrak{m}}
\renewcommand{\AA}{\mathfrak{a}}
\newcommand{\Qpbar}{\overline{\mathbf{Q}}_p}

\renewcommand{\phi}{\varphi}
\renewcommand{\projlim}{\varprojlim}
\renewcommand{\geq}{\geqslant}
\renewcommand{\leq}{\leqslant} 
\newcommand{\NP}{\mathcal{N}_P}
\newcommand{\TP}{\mathcal{T}_P}
\newcommand{\Nm}{\mathrm{N}}
\newcommand{\Tr}{\mathrm{Tr}}
\newcommand{\dS}{\mathrm{S}}
\newcommand{\bigO}{\mathrm{O}}

\newcommand{\calF}{\mathcal{F}}
\newcommand{\Gal}{\mathrm{Gal}}

\newcommand{\End}{\mathrm{End}}
\newcommand{\LT}{\mathrm{LT}}
\newcommand{\Id}{\mathrm{Id}}

\newcommand{\Card}{\mathrm{Card}}
\newcommand{\Fil}{\mathrm{Fil}}

\newcommand{\Col}{\mathrm{Col}}
\newcommand{\val}{\mathrm{val}}

\newcommand{\ab}{\mathrm{ab}}
\newcommand{\unr}{\mathrm{unr}}
\newcommand{\cyc}{\mathrm{cyc}}

\newcommand{\Log}{\mathrm{L}_P}

\newcommand{\at}{\tilde{\mathbf{A}}}
\newcommand{\et}{\widetilde{\mathbf{E}}}
\newcommand{\atplus}{\tilde{\mathbf{A}}^+}
\newcommand{\etplus}{\widetilde{\mathbf{E}}^+}
\newcommand{\btrigplus}{\widetilde{\mathbf{B}}^+_{\mathrm{rig}}}

\newcommand{\bdr}{\mathbf{B}_{\mathrm{dR}}}  
\newcommand{\bcris}{\mathbf{B}_{\mathrm{cris}}}

\newcommand{\dcroc}[1]{[\![ #1 ]\!]}

\newcommand{\ubar}{\overline{u}}

\author{Laurent Berger}
\address{UMPA de l'ENS de Lyon \\
UMR 5669 du CNRS \\ IUF}
\email{laurent.berger@ens-lyon.fr}
\urladdr{perso.ens-lyon.fr/laurent.berger/}
 
\date{June 2015}

\title[Iterated extensions and relative Lubin-Tate groups]{Iterated extensions and \\ relative Lubin-Tate groups}
\alttitle{Extensions it\'er\'ees et groupes de Lubin-Tate relatifs}

\begin{document}

\begin{abstract}
Let $K$ be a finite extension of $\Qp$ with residue field $\FF_q$ and let $P(T) = T^d + a_{d-1}T^{d-1} + \cdots +a_1 T$ where $d$ is a power of $q$ and $a_i \in \MM_K$ for all $i$. Let $u_0$ be a uniformizer of $\OO_K$ and let $\{u_n\}_{n \geq 0}$ be a sequence of elements of $\Qpbar$ such that $P(u_{n+1}) = u_n$ for all $n \geq 0$. Let $K_\infty$ be the field generated over $K$ by all the $u_n$. If $K_\infty / K$ is a Galois extension, then it is abelian, and our main result is that it is generated by the torsion points of a relative Lubin-Tate group (a generalization of the usual Lubin-Tate groups). The proof of this involves generalizing the construction of Coleman power series, constructing some $p$-adic periods in Fontaine's rings, and using local class field theory.
\end{abstract}

\begin{altabstract}
Soit $K$ une extension finie de $\Qp$ de corps r\'esiduel $\FF_q$ et $P(T) = T^d + a_{d-1}T^{d-1} + \cdots +a_1 T$ o\`u $d$ est une puissance de $q$ et $a_i \in \MM_K$ pour tout $i$. Soit $u_0$ une uniformisante de $\OO_K$ et $\{u_n\}_{n \geq 0}$ une suite d'\'el\'ements de $\Qpbar$ telle que $P(u_{n+1}) = u_n$ pour tout $n \geq 0$. Soit $K_\infty$ l'extension de $K$ engendr\'ee par les $u_n$. Si $K_\infty / K$ est Galoisienne, alors elle est ab\'elienne, et notre r\'esultat principal est qu'elle est engendr\'ee par les points de torsion d'un groupe de Lubin-Tate relatif (une g\'en\'eralisation des groupes de Lubin-Tate usuels). Pour prouver cela, nous g\'en\'eralisons la construction des s\'eries de Coleman, construisons des p\'eriodes $p$-adiques dans les anneaux de Fontaine et utilisons la th\'eorie du corps de classes local.
\end{altabstract}



\thanks{This material is partially based upon work supported by the NSF under Grant No.\ 0932078 000 while the author was in residence at the MSRI in Berkeley, California, during the Fall 2014 semester.}

\dedicatory{To Glenn Stevens, on the occasion of his 60th birthday}

\maketitle

\tableofcontents

\setlength{\baselineskip}{18pt}

\section*{Introduction}\label{intro} 

Let $K$ be a field, let $P(T) \in K[T]$ be a polynomial of degree $d \geq 1$, choose $u_0 \in K$ and for $n \geq 0$, let $u_{n+1} \in \overline{K}$ be such that $P(u_{n+1}) = u_n$. The field $K_\infty$ generated over $K$ by all the $u_n$ is called an \emph{iterated extension} of $K$. These iterated extensions and the resulting Galois groups have been studied in various contexts, see for instance \cite{O85}, \cite{S92}, \cite{AHM} and \cite{BJ}.

In this article, we focus on a special situation: $p \neq 2$ is a prime number, $K$ is a finite extension of $\Qp$, with ring of integers $\OO_K$, whose maximal ideal is $\MM_K$ and whose residue field is $k$. Let $d$ be a power of $\Card(k)$, and let $P(T) = T^d + a_{d-1} T^{d-1} + \cdots + a_1 T$ be a monic polynomial of degree $d$ with $a_i \in \MM_K$ for $1 \leq i \leq d-1$. Let $u_0$ be a uniformizer of $\OO_K$ and define a sequence $\{ u_n \}_{n \geq 0}$ by letting $u_{n+1}$ be a root of $P(T) =  u_n$. Let $K_n = K(u_n)$ and $K_\infty=\cup_{n \geq 1} K_n$. This iterated extension is called a \emph{Frobenius-iterate extension}, after \cite{CD14} (whose definition is a bit more general than ours). The question that we consider in this article is: which \emph{Galois} extensions $K_\infty/K$ are Frobenius iterate? 

This question is inspired by the observation, made in remark 7.16 of \cite{CD14}, that it follows from the main results of ibid.\  and \cite{LTFN} that: if $K_\infty/K$ is Frobenius-iterate and Galois, then it is necessarily abelian. Here, we prove a much more precise result. 

First, let us recall that in \cite{DSLT}, de Shalit gives a generalization of the construction of Lubin-Tate formal groups (for which see \cite{LT}). A \emph{relative Lubin-Tate group} is a formal group $\dS$ that is attached to an unramified extension $E/F$ and to an element $\alpha$ of $F$ of valuation $[E:F]$. The extension $E_\infty^{\dS}/F$ generated over $F$ by the torsion points of this formal group is the subextension of $F^{\ab}$ cut out via local class field theory by the subgroup of $F^\times$ generated by $\alpha$. If $E=F$, we recover the classical Lubin-Tate groups.

\begin{enonce*}{Theorem A}
Let $K$ be a finite Galois extension of $\Qp$, and let $K_\infty/K$ be a Frobenius-iterate extension. If $K_\infty/K$ is Galois, then there exists a subfield $F$ of $K$, and a relative Lubin-Tate group $\dS$, relative to the extension $F^{\unr} \cap K$ of $F$, such that if $K_\infty^{\dS}$ denotes the extension of $K$ generated by the torsion points of $\dS$, then $K_\infty \subset K_\infty^{\dS}$ and $K_\infty^{\dS} / K_\infty$ is a finite extension.
\end{enonce*}

This is theorem \ref{dslt}. Conversely, it is easy to see that the extension coming from a relative Lubin-Tate group is Frobenius-iterate after the first layer (see example \ref{expolit}). The proof of theorem A is quite indirect. We start with the observation that if $K_\infty / K$ is a Frobenius-iterate extension, that is not necessarily Galois, then we can generalize the construction of Coleman's power series (see \cite{C79}). Let $\projlim \OO_{K_n}$ denote the set of sequences $\{x_n\}_{n \geq 0}$ with $x_n \in \OO_{K_n}$ and such that $\Nm_{K_{n+1}/K_n}(x_{n+1}) = x_n$ for all $n \geq 0$.

\begin{enonce*}{Theorem B}
We have $\{u_n\}_{n \geq 0} \in \projlim \OO_{K_n}$ and if $\{x_n\}_{n \geq 0} \in \projlim \OO_{K_n}$, then there exists a unique power series $\Col_x(T) \in \OO_K \dcroc{T}$ such that $x_n = \Col_x(u_n)$ for all $n \geq 0$.
\end{enonce*}

Suppose now that $K_\infty/K$ is Galois, and let $\Gamma=\Gal(K_\infty/K)$. The results of \cite{CD14} and \cite{LTFN} imply that $K_\infty/K$ is abelian, so that $K_n/K$ is Galois for all $n \geq 1$. If $g \in \Gamma$, then $\{g(u_n)\}_{n \geq 0} \in \projlim \OO_{K_n}$, so that by theorem B, we get a power series $\Col_g(T) \in \OO_K \dcroc{T}$ such that $g(u_n) = \Col_g(u_n)$ for all $n \geq 0$. Let $\etplus = \projlim_{x \mapsto x^d} \OO_{\Cp} / p$, let $K_0 = \Qp^{\unr} \cap K$ and let $\atplus = \OO_K \otimes_{\OO_{K_0}} W(\etplus)$ be Fontaine's rings of periods (see \cite{FP}). The element $\{u_n\}_{n \geq 0}$ gives rise to an element $\ubar \in \etplus$.

\begin{enonce*}{Theorem C}
There exists $u \in \atplus$ whose image in $\etplus$ is $\ubar$, and such that $\phi_d(u) = P(u)$. We have $g(u)= \Col_g(u)$ if $g \in \Gamma$.
\end{enonce*}

The power series $\Col_g(T)$ satisfies the functional equation $\Col_g \circ P(T) = P \circ \Col_g(T)$. The study of $p$-adic power series that commute under composition was taken up by Lubin in \cite{L94}. In \S 6 of ibid.,  Lubin writes that ``experimental evidence seems to suggest that for an invertible series to commute with a noninvertible series, there must be a formal group somehow in the background''. There are a number of results in this direction, see for instance \cite{LMS}, \cite{SS13} and \cite{JS14}. In our setting, the series $\{ \Col_g(T) \}_{g \in \Gamma}$ commute with $P(T)$ and theorem A says that indeed, there is a formal group that accounts for this. Let us now return to the proof of theorem A. We first show that $P'(T) \neq 0$. It is then proved in \S 1 of \cite{L94} that given such a $P(T)$, a power series $\Col_g(T)$ that commutes with $P(T)$ is determined by $\Col_g'(0)$. If we let $\eta(g) = \Col_g'(0)$, we get the following: the map $\eta : \Gamma \to \OO_K^\times$ is an injective character.

In order to finish the proof of theorem A, we use some $p$-adic Hodge theory. Let $\Log(T) \in K\dcroc{T}$ be the \emph{logarithm} attached to $P(T)$ and constructed in \cite{L94}; it converges on the open unit disk, and satisfies $\Log \circ P(T) = P'(0) \cdot \Log(T)$ as well as  $\Log \circ \Col_g(T) = \eta(g) \cdot \Log(T)$ for $g \in \Gamma$. In particular, we can consider $\Log(u)$ as an element of the ring $\bcris^+$ (see \cite{FP} for the rings of periods $\bcris^+$ and $\bdr$), which satisfies $g(\Log(u)) = \eta(g) \cdot \Log(u)$. More generally, if $\tau \in \Gal(K/\Qp)$, then we can twist $u$ by $\tau$ to get some elements $u_\tau \in \atplus$ and $\Log^\tau(u_\tau) \in \bcris^+$, satisfying $g(\Log^\tau(u_\tau)) = \tau(\eta(g)) \cdot \Log^\tau(u_\tau)$. The elements $\{\Log^\tau(u_\tau)\}_\tau$ are crystalline periods for the representation arising from $\eta$. Our main technical result concerning these periods is that the set of $\tau \in \Gal(K/\Qp)$ such that $\Log^\tau(u_\tau) \in \Fil^1 \bdr$ is a subgroup of $\Gal(K/\Qp)$, and therefore cuts out a subfield $F$ of $K$. This allows us to prove the following.

\begin{enonce*}{Theorem D}
There exists a subfield $F$ of $K$, a Lubin-Tate character $\chi_K$ attached to a uniformizer of $K$, and an integer $r \geq 1$, such that $\eta = \Nm_{K/F}(\chi_K)^r$.
\end{enonce*}

Theorem A follows from theorem D by local class field theory: the extensions of $K$ corresponding to $\Nm_{K/F}(\chi_K)$ are precisely those that come from relative Lubin-Tate groups. At the end of \S \ref{lcft}, we give an example for which $r=2$. In this example, the Coleman power series $p$-adically  interpolate Chebyshev polynomials.

\section{Relative Lubin-Tate groups}
\label{dsltsec}

We recall de Shalit's construction (see \cite{DSLT}) of a family of formal groups that generalize Lubin-Tate groups. Let $F$ be a finite extension of $\Qp$, with ring of integers $\OO_F$ and residue field $k_F$ of cardinality $q$. Take $h \geq 1$ and let $E$ be the unramified extension of $F$ of degree $h$. Let $\phi_q : E \to E$ denote the Frobenius map that lifts $[x \mapsto x^q]$. If $f(T) = \sum_{i \geq 0} f_i T^i \in E\dcroc{T}$, let $f^{\phi_q}(T) = \sum_{i \geq 0} \phi_q(f_i) T^i$. 

If $\alpha \in \OO_F$ is such that $\val_F(\alpha) = h$, let $\calF_\alpha$ be the set of power series $f(T) \in \OO_E \dcroc{T}$ such that $f(T) = \pi T + \bigO(T^2)$ with $\Nm_{E/F}(\pi) = \alpha$ and such that $f(T) \equiv T^q \bmod{\MM_E \dcroc{T}}$. The set $\calF_\alpha$ is nonempty, since $\Nm_{E/F}(E^\times)$ is the set of elements of $F^\times$ whose valuation is in $h \cdot \ZZ$. If $\Nm_{E/F}(\pi) = \alpha$, one can take $f(T)=\pi T+T^q$. The following theorem summarizes some of the results of \cite{DSLT} (see also \S IV of \cite{IK}).

\begin{theo}\label{dsltmain}
If $f(T) \in \calF_\alpha$, then 
\begin{enumerate}
\item there is a unique formal group law $\dS(X,Y) \in \OO_E \dcroc{X,Y}$ such that $\dS^{\phi_q} \circ f = f \circ \dS$, and the isomorphism class of $\dS$ depends only on $\alpha$;
\item for all $a \in \OO_F$, there exists a unique power series $[a](T) \in \OO_E \dcroc{T}$ such that $[a](T) = aT + \bigO(T^2)$ and $[a](T) \in \End(\dS)$.
\end{enumerate}
Let $x_0 = 0$ and for $m \geq 0$, let $x_m \in \Qpbar$ be such that $f^{\phi_q^m}(x_{m+1}) = x_m$ (with $x_1 \neq 0$). Let $E_m = E(x_m)$ and let $E_\infty^{\dS} = \cup_{m \geq 1} E_m$. 
\begin{enumerate}
\item The fields $E_m$ depend only on $\alpha$, and not on the choice of $f(T) \in \calF_\alpha$;
\item The extension $E_m/E$ is Galois, and its Galois group is isomorphic to $(\OO_F/\MM_F^m)^\times$;
\item $E_\infty^{\dS} \subset F^{\ab}$ and $E_\infty^{\dS}$ is the subfield of $F^{\ab}$ cut out by $\langle \alpha \rangle \subset F^\times$ via local class field theory.
\end{enumerate}
\end{theo}

\begin{rema}\label{dshone}
If $h=1$, then we recover the usual Lubin-Tate formal groups of \cite{LT}.
\end{rema}

\section{Frobenius-iterate extensions}
\label{polit} 

Let $p \neq 2$ be a prime number, let $K$ be a finite extension of $\Qp$, with ring of integers $\OO_K$, whose maximal ideal is $\MM_K$ and whose residue field is $k$. Let $q=\Card(k)$, and let $\pi$ denote a uniformizer of $\OO_K$. Let $d$ be a power of $q$, and let $P(T) = T^d + a_{d-1} T^{d-1} + \cdots + a_1 T$ be a monic polynomial of degree $d$ with $a_i \in \MM_K$ for $1 \leq i \leq d-1$.

Let $u_0$ be a uniformizer of $\OO_K$ and define a sequence $\{ u_n \}_{n \geq 0}$ by letting $u_{n+1}$ be a root of $P(T) =  u_n$. Let $K_n = K(u_n)$.

\begin{lemm}\label{pinunif}
The extension $K_n/K$ is totally ramified of degree $d^n$, $u_n$ is a uniformizer of $\OO_{K_n}$ and $\Nm_{K_{n+1}/K_n} ( u_{n+1} ) = u_n$.
\end{lemm}

\begin{proof}
The first two assertions follow immediately from the theory of Newton polygons, and the last one from the fact that $P(T) - u_n$ is the minimal polynomial of $u_{n+1}$ over $K_n$, as well as the fact that $d$ is odd since $p \neq 2$. 
\end{proof}

Let $K_\infty = \cup_{n \geq 1} K_n$. This is a totally ramified infinite and pro-$p$ extension of $K$.

\begin{defi}\label{defpolit}
We say that an extension $K_\infty/K$ is \emph{$\varphi$-iterate} if it is of the form above.
\end{defi}

This definition is inspired by the similar one that is given in definition 1.1 of \cite{CD14}. We require $P(T)$ to be a monic polynomial, instead of a more general power series as in ibid., in order to control the norm of $u_n$ and to ensure the good behavior of $K_{n+1}/K_n$.

\begin{exem}\label{expolit}
(i) If $P(T)=T^q$, then $K_\infty/K$ is a $\phi$-iterate extension, which is the Kummer extension of $K$ corresponding to $\pi$.

(ii) Let $\LT$ be a Lubin-Tate formal $\OO_K$-module attached to $\pi$, and $K_n=K(\LT[\pi^n])$. The extension $K_\infty/K_1$ is $\varphi$-iterate with $P(T) = [\pi](T)$. 

(iii) More generally, let $\dS$ be a relative Lubin-Tate group, relative to an extension $E/F$ and $\alpha \in F$ as in \S \ref{dsltsec}. The extension $E_\infty^{\dS}/E_1$ is $\varphi$-iterate with $P(T) = [\alpha](T)$.
\end{exem}

\begin{proof}
Item (ii) follows from applying (iii) with $K=E=F$, and we now prove (iii). We use the notation of theorem \ref{dsltmain}. Since the isomorphism class of $\dS$ and the extension $E_\infty^{\dS}/E$ only depend on $\alpha$, we can take $f(T) = \pi T + T^q$ where $\Nm_{E/F}(\pi) = \alpha$. Let $P(T) = f^{\phi_q^{h-1}} \circ \cdots \circ f^{\phi_q} \circ f(T)  \in \OO_E[T]$, so that $P(T) = [\alpha](T)$. The extension $E_{hm+1}$ is generated by $x_{hm+1}$ over $E_1$, and we have $P(x_{hm+1}) = x_{(h-1)m+1}$.  The claim therefore follows from taking $u_m = x_{hm+1}$ for $m \geq 0$, and observing that since $\pi + u_0^{q-1} = 0$, $u_0$ is a uniformizer of $\OO_{E_1}$.
\end{proof}

\section{Coleman power series}
\label{colpow}

Let us write $\projlim \OO_{K_n}$ for the set of sequences $\{x_n\}_{n \geq 0}$ such that $x_n \in \OO_{K_n}$ and such that $\Nm_{K_{n+1}/K_n}(x_{n+1})=x_n$ for $n \geq 0$. By lemma \ref{pinunif}, the sequence $\{ u_n \}_{n \geq 0}$ belongs to $\projlim \OO_{K_n}$. The goal of this {\S} is to show the following theorem (theorem B).

\begin{theo}\label{colexist}
If $\{ x_n \}_{n \geq 0} \in \projlim \OO_{K_n}$, then there exists a uniquely determined power series $\Col_x(T) \in \OO_K \dcroc{T}$ such that $x_n = \Col_x(u_n)$ for all $n \geq 0$.
\end{theo}

Our proof follows the one that is given in \S 13 of \cite{W97}. The unicity is a consequence of the following well-known general principle.

\begin{prop}
\label{colunik}
If $f(T) \in \OO_K \dcroc{T}$ is nonzero, then $f(T)$ has only finitely many zeroes in the open unit disk.
\end{prop}

In order to prove the existence part of theorem \ref{colexist}, we start by generalizing Coleman's norm map (see \cite{C79} for the original construction, and \S 2.3 of \cite{F90} for the generalization that we use). The ring $\OO_K\dcroc{T}$ is a free $\OO_K\dcroc{P(T)}$-module of rank $d$. If $f(T) \in \OO_K\dcroc{T}$, let $\NP(f)(T) \in \OO_K\dcroc{T}$ be defined by the requirement that $\NP(f)(P(T)) = \Nm_{\OO_K\dcroc{T} / \OO_K\dcroc{P(T)}} (f(T))$. For example, $\NP(T)=T$ since $d$ is odd.

\begin{prop}\label{pronmp}
The map $\NP$ has the following properties.
\begin{enumerate}
\item If $f(T) \in \OO_K\dcroc{T}$, then $\NP(f)(u_n) = \Nm_{K_{n+1}/K_n} (f(u_{n+1}))$; 
\item If $k \geq 1$ and $f(T) \in 1+ \pi^k \OO_K\dcroc{T}$, then $\NP(f)(T) \in 1+ \pi^{k+1} \OO_K\dcroc{T}$;
\item If $f(T) \in \OO_K \dcroc{T}$, then $\NP(f)(T) \equiv f(T) \bmod{\pi}$;
\item If $f(T) \in \OO_K \dcroc{T}^\times$, and $k,m \geq 0$, then $\NP^{m+k}(f) \equiv \NP^k(f) \bmod{\pi^{k+1}}$.
\end{enumerate}
\end{prop}

\begin{proof}
The determinant of the multiplication-by-$f(T)$ map on the $\OO_K\dcroc{P(T)}$-module $\OO_K\dcroc{T}$ is $\NP(f)(P(T))$. By evaluating at $T=u_{n+1}$, we find that the determinant of the multiplication-by-$f(u_{n+1})$ map on the $\OO_{K_n}$-module $\OO_{K_{n+1}}$ is $\NP(f)(u_n)$, so that $\NP(f)(u_n) = \Nm_{K_{n+1}/K_n} (f(u_{n+1}))$.

We now prove (2). If $f(T) \in \OO_K\dcroc{T}$, let $\TP(f)(T) \in \OO_K\dcroc{T}$ be the trace map defined by $\TP(f)(P(T)) = \Tr_{\OO_K\dcroc{T} / \OO_K\dcroc{P(T)}} (f(T))$. A straightforward calculation shows that if $h(T) \in \OO_K \dcroc{T}$, then $\TP(h)(T) \in \pi \cdot \OO_K \dcroc{T}$. If $f(T) = 1 +\pi^k h(T)$, then $\NP(f)(T) \equiv 1 + \pi^k \TP(h)(T) \bmod{\pi^{k+1}}$, so that $\NP(f)(T) \in 1+ \pi^{k+1} \OO_K\dcroc{T}.$

Item (3) follows from a straightfoward calculation in $k\dcroc{T}$ using the fact that $P(T)=T^d$ in $k\dcroc{T}$. Finally, let us prove (4). If $f(T) \in \OO_K \dcroc{T}^\times$, then $\NP(f)/f \equiv 1 \bmod{\pi}$ by (3), so that $\NP^m(f)/f \equiv 1 \bmod{\pi}$ as well. Item (2) now implies that $\NP^{m+k}(f) \equiv \NP^k(f) \bmod{\pi^{k+1}}$.
\end{proof}

\begin{proof}[of theorem \ref{colexist}]
The power series $\Col_x(T)$ is unique by lemma \ref{colunik}, and we now show its existence. 
If $x_n$ is not a unit of $\OO_{K_n}$, then there exists $e \geq 1$ such that $x_n=u_n^e x_n^*$ where $x_n^* \in \OO_{K_n}^\times$ for all $n$, and then $\Col_x(T) = T^e \cdot \Col_{x^*}(T)$. We can therefore assume that $x_n$ is a unit of $\OO_{K_n}$. For all $j \geq 1$, we have  $\OO_{K_j} = \OO_K[u_j]$, so that there exists $g_j(T) \in \OO_K[T]$ such that $x_j = g_j(u_j)$. Let $f_j(T) = \NP^j(g_{2j})$. By proposition \ref{pronmp}, we have $x_n \equiv f_j(u_n) \bmod{\pi^{j+1}}$ for all $n \leq j$. The space $\OO_K \dcroc{T}$ is compact; let $f(T)$ be a limit point of $\{ f_j \}_{j \geq 1}$. We have $x_n = f(u_n)$ for all $n$ by continuity, so that we can take $\Col_x(T) = f(T)$.
\end{proof}

\begin{rema}\label{colnpinv}
We have $\NP(\Col_x)(T) = \Col_x(T)$.
\end{rema}

\begin{proof}
The power series $\NP(\Col_x)(T) - \Col_x(T)$ is zero at $T=u_n$ for all $n \geq 0$ by proposition \ref{pronmp}, so that $\NP(\Col_x)(T) = \Col_x(T)$ by lemma \ref{colunik}.
\end{proof}

\section{Lifting the field of norms}
\label{cdn}

In this {\S}, we assume that $K_\infty/K$ is a Galois extension, and let $\Gamma = \Gal(K_\infty / K)$. We recall some results of \cite{CD14} and \cite{LTFN}, and give a more precise formulation of some of them in our specific situation.

\begin{prop}\label{kngal}
If $K_\infty/K$ is Galois, then $K_n/K$ is Galois for all $n \geq 1$.
\end{prop}

\begin{proof}
It follows from the main results of \cite{CD14} and of \cite{LTFN} (see remark 7.16 of \cite{CD14}) that if $K_\infty/K$ is a $\phi$-iterate extension that is Galois, then it is abelian. This implies the proposition (it would be more satisfying to find a direct proof).
\end{proof}

If $g \in \Gamma$, proposition \ref{kngal} and theorem \ref{colexist} imply that there is a unique power series $\Col_g(T) \in \OO_K \dcroc{T}$ such that $g(u_n) = \Col_g(u_n)$ for all $n \geq 0$. In the sequel, we need some ramification-theoretic properties of $K_\infty/K$. They are summarized in the theorem below.

\begin{theo}\label{normpow}
There exists a constant $c=c(K_\infty/K) > 0$ such that for any $E \subset F$, finite extensions of $K$ contained in $K_\infty$,  and $x \in \OO_F$, we have
\[ \val_K \left( \frac{\Nm_{F/E}(x)}{x^{[F:E]}}-1 \right) \geq c. \]
\end{theo}

\begin{proof}
By the main result of \cite{CDAPF}, the extension $K_\infty/K$ is \emph{strictly APF}, so that if we denote by $c(K_\infty/K)$ the constant defined in 1.2.1 of \cite{W83}, then $c(K_\infty/K) > 0$. By 4.2.2.1 of ibid., we have
\[ \val_E \left( \frac{\Nm_{F/E}(x)}{x^{[F:E]}}-1 \right) \geq c(F/E), \]
By 1.2.3 of ibid., $c(F/E) \geq c(K_\infty/E)$ and (see for instance the proof of 4.5 of \cite{CD14} or page 83 of \cite{W83}) $c(K_\infty/E) \geq c(K_\infty/K) \cdot [E:K]$. This proves the theorem.
\end{proof}

Let $c$ be the constant afforded by theorem \ref{normpow}. We can always assume that $c \leq \val_K(p)/(p-1)$. If $E$ is some subfield of $\Cp$, let $\AA_E^c$ denote the set of elements $x$ of $E$ such that $\val_K(x) \geq c$. Let $\et^+= \projlim_{x \mapsto x^d} \OO_{\Cp} / \AA^c_{\Cp}$. The sequence $\{ u_n \}_{n \geq 0}$ gives rise to an element $\ubar \in \etplus$. Recall that by \S 2.1 and \S 4.2 of \cite{W83}, there is an embedding $\iota : \projlim \OO_{K_n} \to \etplus$, which is an isomorphism onto $\projlim_{x \mapsto x^d} \OO_{K_n} / \AA^c_{K_n}$, which is also isomorphic to $k \dcroc{\ubar}$. Let $K_0= \Qp^{\unr} \cap K$ and $\at^+=\OO_K \otimes_{\OO_{K_0}} W(\et^+)$. Recall (see \cite{FP}) that we have a map $\theta : \atplus \to \OO_{\Cp}$. If $x \in \atplus$ and $\overline{x} = (x_n)_{n \geq 0}$ in $\etplus$, then $\theta \circ \phi_d^{-n}(x) = x_n$ in $\OO_{\Cp} /  \AA^c_{\Cp}$.

\begin{theo}\label{liftubar}
There exists a unique $u \in \atplus$ whose image in $\etplus$ is $\ubar$, and such that $\phi_d(u) = P(u)$. Moreover: 
\begin{enumerate}
\item[(i)] If $n \geq 0$, then $\theta \circ \phi_d^{-n}(u) = u_n$;
\item[(ii)] $\OO_K \dcroc{u} = \{ x \in \atplus$, $\theta \circ \phi_d^{-n}(x) \in \OO_{K_n}$ for all $n \geq 1\}$;
\item[(iii)] $g(u) = \Col_g(u)$ if $g \in \Gamma$.
\end{enumerate}
\end{theo}

\begin{proof}
The existence of $u$ and item (i) are proved in lemma 9.3 of \cite{CEV}, where it is shown that $u = \lim_{n \to +\infty} P^{\circ n}(\phi_d^{-n}([\ubar]))$.

Let $R =  \{ x \in \atplus$ such that $\theta \circ \phi_d^{-n}(x) \in \OO_{K_n}$ for all $n \geq 1\}$. If $x \in R$, then its image in $\etplus$ lies in $\projlim_{x \mapsto x^d} \OO_{K_n} / \AA^c_{K_n} = k \dcroc{\ubar}$. We have $u \in R$ by proposition \ref{liftubar}, so that the map $R / \pi R \to k \dcroc{\ubar}$ is surjective. We then have $R = \OO_K \dcroc{u}$, since $R$ is separated and complete for the $\pi$-adic topology, which proves (ii).

The ring $\OO_K \dcroc{u}$ is stable under the action of $G_K$ by (ii). If $g \in \Gamma$, there exists $F_g(T) \in \OO_K \dcroc{T}$ such that $g(u) = F_g(u)$. We have $g(u_n) = g ( \theta \circ \phi_d^{-n}(u)) = \theta \circ \phi_d^{-n} (F_g(u)) = F_g(u_n)$ by (i), so that  $g(u_n) = F_g(u_n)$ for all $n$. This implies that $F_g(T) = \Col_g(T)$.
\end{proof}

\begin{rema}\label{liftfn}
In the terminology of \cite{W83}, $\projlim \OO_{K_n}$ is the ring of integers of the field of norms $X(K_\infty)$ of the extension $K_\infty/K$, and theorem \ref{liftubar} shows that we can lift $X(K_\infty)$ to characteristic zero, along with the Frobenius map $\phi_d$ and the action of $\Gamma$.
\end{rema}

If $g \in \Gamma$, then $\Col_g \circ P(T) = P \circ \Col_g(T)$ since the two series have the same value at $u_n$ for all $n \geq 1$. Let $\eta(g) = \Col_g'(0)$, so that $g \mapsto \eta(g)$ is a character $\eta : \Gamma \to \OO_K^\times$

\begin{prop}\label{comnoz}
If $F(T) \in T \cdot \OO_K \dcroc{T}$ is such that $F'(0) \in 1+p \OO_K$, and if $A(T) \in T \cdot \OO_K \dcroc{T}$ vanishes at order $k \geq 2$ at $0$, and satisifes $A \circ F(T) = F \circ A(T)$, then $F(T) = T$.
\end{prop}

\begin{proof}
Write $F(T) = f_1 T + \bigO(T^2)$, and $A(T) = a_k T^k + \bigO(T^{k+1})$ with $a_k \neq 0$. The equation $F \circ A(T) = A \circ F(T)$ implies that $f_1 a_k = a_k f_1^k$ so that if $k \neq 1$, then $f_1^{k-1}=1$. Since $f_1 \in 1+p\OO_K$ and $p \neq 2$, this implies that $f_1=1$. If $F(T) \neq T$, we can write $F(T) = T + T^i h(T)$ for some $i \geq 2$ with $h(0) \neq 0$. The equation $F \circ A(T) = A \circ F(T)$ and the equality $A(T + T^i h(T))=\sum_{j \geq 0} (T^i h(T))^j A^{(j)}(T)/j!$ imply that $A(T) + A(T)^i h(A(T))  = A(T) + T^i h(T) A'(T) + \bigO (T^{2i+k-2})$, so that $A(T)^i h(A(T))  = T^i h(T) A'(T) + \bigO (T^{2i+k-2})$. The term of lowest degree of the LHS is of degree $ki$, while on the RHS it is of degree $i+k-1$. We therefore have $ki=i+k-1$, so that $(k-1)(i-1)=0$ and hence $k=1$.
\end{proof}

\begin{coro}\label{dernonul}
We have $P'(0) \neq 0$.
\end{coro}

\begin{proof}
This follows from proposition \ref{comnoz}, since $\Col_g'(0) \in 1+p\OO_K$ if $g$ is close enough to $1$, and $\Col_g(T) = T$ if and only if $g=1$ (compare with lemma 4.5 of \cite{LTFN}).
\end{proof}

\begin{coro}\label{galisab}
The character $\eta : \Gamma \to \OO_K^\times$ is injective.
\end{coro}

\begin{proof}
This follows from proposition 1.1 of \cite{L94}, which says that if $P'(0) \in \MM_K \setminus \{ 0 \}$, then a power series $F(T) \in T \cdot \OO_K\dcroc{T}$ that commutes with $P(T)$ is determined by $F'(0)$. This implies that $\Col_g(T)$ is determined by $\eta(g)$, and then $g$ itself is determined by $\Col_g(T)$, since $g(u_n) = \Col_g(u_n)$ for all $n$.
\end{proof}

We therefore have a character $\eta : \Gal(\Qpbar/K) \to \OO_K^\times$, such that $K_\infty = \Qpbar^{\ker \eta}$.

\section{$p$-adic Hodge theory}
\label{chareta}

We now assume that $K/\Qp$ is Galois (for simplicity), and we keep assuming that $K_\infty/K$ is Galois. We use the element $u$ above, and Lubin's logarithm (proposition \ref{lublog} below), to construct crystalline periods for $\eta$.

\begin{prop}\label{lublog}
There exists a power series $\Log(T) \in K\dcroc{T}$ that is holomorphic on the open unit disk, and satisfies 
\begin{enumerate}
\item $\Log(T) = T + \bigO(T^2)$;
\item $\Log \circ P(T) = P'(0) \cdot \Log(T)$;
\item $\Log \circ \Col_g(T) = \eta(g) \cdot \Log(T)$ if $g \in \Gamma$.
\end{enumerate}
If we write $P(T)=T \cdot Q(T)$, then 
\[ \Log(T) = \lim_{n \to +\infty} \frac{P^{\circ n}(T)} {P'(0)^n} = T \cdot \prod_{n \geq 0} \frac{Q(P^{\circ n}(T))}{Q(0)}. \]
\end{prop}

\begin{proof}
See propositions 1.2, 2.2 and 1.3 of \cite{L94}.
\end{proof}

Let $\btrigplus$ denote the Fr\'echet completion of $\atplus [1/\pi]$, so that our $\btrigplus$ is $K \otimes_{K_0}$ the ``usual'' $\btrigplus$ (for which see \cite{LB2}). If $u \in \atplus$ is the element afforded by proposition \ref{liftubar}, then $\Log(u)$ converges in $\btrigplus$. We have $g(\Log(u)) = \eta(g) \cdot \Log(u)$ by proposition \ref{lublog}. If $\tau \in \Gal(K/\Qp)$, then let $n(\tau)$ be some $n \in \ZZ$ such that $\tau = \phi^n$ on $k_K$, and let $u_\tau = (\tau \otimes \phi^{n(\tau)}) (u) \in \atplus$. 

If $F(T) = \sum_{i \geq 0} f_i T^i \in K \dcroc{T}$, let $F^\tau(T) = \sum_{i \geq 0} \tau(f_i) T^i$. We have $g(\Log^\tau(u_\tau)) = \tau(\eta(g)) \cdot \Log^\tau(u_\tau)$ in $\btrigplus$. This implies the following result, which is a slight improvement of theorem 4.1 of \cite{LTFN}.

\begin{prop}\label{cisdrp}
The character $\eta : \Gamma \to \OO_K^\times$ is crystalline, with weights in $\ZZ_{\geq 0}$.
\end{prop}

\begin{proof}
The fact that $g(\Log^\tau(u_\tau)) = \tau(\eta(g)) \cdot \Log^\tau(u_\tau)$ for all $\tau \in \Gal(K/\Qp)$ implies that $\eta$ gives rise to a $K \otimes_{K_0} \bcris$-admissible representation. If $V$ is any $p$-adic representation of $G_K$, then
\[ \left((K \otimes_{K_0} \bcris) \otimes_{\Qp} V \right)^{G_K} = K \otimes_{K_0} (\bcris \otimes_{\Qp} V)^{G_K}. \]
This implies that a $K \otimes_{K_0} \bcris$-admissible representation is crystalline. The weights of $\eta$ are $\geq 0$ because $\Log^\tau(u_\tau) \in \bdr^+$ for all $\tau$.
\end{proof}

\begin{lemm}\label{thetutau}
We have $\theta \circ \phi_d^{-n}(u_\tau) = \lim_{k \to +\infty} (P^\tau)^{\circ k}(u_{n+k}^{p^{n(\tau)}})$.
\end{lemm}

\begin{proof}
The element $u_\tau \in \atplus$ has the property that its image in $\etplus$ is $\phi^{n(\tau)}(\ubar) = \ubar^{p^{n(\tau)}}$, and that $\phi_d(u_\tau) = P^\tau(u_\tau)$. The lemma then follows from lemma 9.3 of \cite{CEV}.
\end{proof}

For simplicity, write $u_\tau^n = \theta \circ \phi_d^{-n}(u_\tau)$ and $u_\tau^{n,k} = (P^\tau)^{\circ k}(u_{n+k}^{p^{n(\tau)}})$.

\begin{lemm}\label{unifconv}
If $M > 0$, there exists $j \geq 0$ such that $\val_K(u_\tau^n - u_\tau^{n,j}) \geq M$ for $n \geq 1$.
\end{lemm}

\begin{proof}
If $c$ is the constant coming from theorem \ref{normpow}, then $\val_K(u_\tau^n - u_\tau^{n,0}) \geq c$ for all $n \geq 1$. We prove the lemma by inductively constructing a sequence $\{c_j\}_{j \geq 0}$ such that $\val_K(u_\tau^n - u_\tau^{n,j}) \geq c_j$ for all $n \geq 1$, and such that $c_j \geq M$ for $j \gg 0$. Let $c_0=c$ and suppose that for some $j$, we have $\val_K(u_\tau^n - u_\tau^{n,j}) \geq c_j$ for all $n \geq 1$. We then have 
\[ \val_K(u_\tau^n - u_\tau^{n,j+1}) = \val_K \left( P^\tau(u_\tau^{n+1}) - P^\tau(u_\tau^{n+1,j})\right). \]

If $R(T) \in \OO_K[T]$ and $x,y \in \OO_{\Qpbar}$, then $R(x)-R(y) = (x-y) R'(y) + (x-y)^2 S(x,y)$ with $S(T,U) \in \OO_K [T,U]$. This, and the fact that $P'(T) \in \MM_K[T]$, implies that we can take $c_{j+1} = \min(c_j+1, 2c_j)$. The lemma follows.
\end{proof}

We now recall a result from \cite{L94}. If $f(T) \in T \cdot \OO_K \dcroc{T}$ is such that $f'(0) \in \MM_K \setminus \{ 0 \}$, let $\Lambda(f)$ be the set of the roots of all iterates of $f$. If $u(T) \in T \cdot \OO_K \dcroc{T}$ is such that $u'(0) \in \OO_K^\times$ and $u'(0)$ is not a root of $1$, let $\Lambda(u)$ be the set of the fixed points of all iterates of $u$.

\begin{lemm}\label{lublamb}
If $f$ and $u$ are as above, and if $u \circ f = f \circ u$, then $\Lambda(f)=\Lambda(u)$.
\end{lemm}

\begin{proof}
This is proposition 3.2 of \cite{L94}.
\end{proof}

For each $\tau \in \Gal(K/\Qp)$, let $r_\tau$ be the weight of $\eta$ at $\tau$. 

\begin{prop}\label{algfil}
If $\tau \in \Gal(K/\Qp)$, then the following are equivalent.
\begin{enumerate}
\item[(i)] $r_\tau \geq 1$;
\item[(ii)]  $\Log^\tau(u_\tau) \in \Fil^1 \bdr$;
\item[(iii)]  $\theta(u_\tau) \in \Qpbar$;
\item[(iv)]  $\theta(u_\tau) \in \Lambda(P^\tau)$;
\item[(v)]  $u_\tau \in \cup_{j \geq 0} \ \phi_d^{-j}(\OO_K \dcroc{u})$.
\end{enumerate}
\end{prop}

\begin{proof}
The equivalence between (i) and (ii) is immediate. We now prove that (ii) implies (iii). If $\Log^\tau(u_\tau) \in \Fil^1 \bdr$, then $\Log^\tau(\theta(u_\tau)) = 0$ so that $\theta(u_\tau) \in \Qpbar$ since it is a root of a convergent power series with coefficients in $K$. We next prove that (iii) implies (iv) (it is clear that (iv) implies (iii)). If $x= \theta(u_\tau)$ then $g(x) = \Col_g^\tau(x)$. If $x \in \Qpbar$ and if $g$ is close enough to $1$, then $g(x)=x$ so that $x \in \Lambda(\Col_g^\tau)$, and then $x \in \Lambda(P^\tau)$ by lemma \ref{lublamb}. Let us prove that (iv) implies (ii). If there exists $n \geq 0$ such that $(P^{\tau})^{\circ n}(\theta(u_\tau)) = 0$, then $(P^{\tau})^{\circ n}(u_\tau) \in \Fil^1 \bdr$ so that $\Log^\tau(u_\tau) \in \Fil^1 \bdr$ as well by proposition \ref{lublog}. Conditions (i), (ii), (iii) and (iv) are therefore equivalent. Condition (v) implies (iii) by using theorem \ref{liftubar} as well as the fact that $\phi_d(u_\tau) = P^\tau(u_\tau)$. 

It remains to prove that (iii) implies (v). Recall that $u_\tau^n = \theta \circ \phi_d^{-n}(u_\tau)$. It is enough to show that there exists $j \geq 0$ such that $u_\tau^n \in \OO_{K_{n+j}}$ for all $n$,  since by theorem \ref{liftubar}, this implies that $u_\tau \in \phi_d^{-j}(\OO_K \dcroc{u})$. Recall that $u_\tau^{n,k} = (P^\tau)^{\circ k}(u_{n+k}^{p^{n(\tau)}})$. Take $M \geq 1 + \val_K((P^\tau)'(u_\tau^n))$ for all $n \gg 0$. By lemma \ref{unifconv}, there exists $j \geq 0$ such that $\val_K(u_\tau^n - u_\tau^{n,j}) \geq M$ for all $n \geq 1$.  The element $u_\tau^n$ is a root of $P^\tau(T) = u_\tau^{n-1}$, and therefore $u_\tau^n- u_\tau^{n,j}$ is a root of $P^\tau(u_\tau^{n,j}+T) - u_\tau^{n-1}$. If $u_\tau^{n-1} \in \OO_{K_{n+j-1}}$, then the polynomial $R_n(T) = P^\tau(u_\tau^{n,j}+T) - u_\tau^{n-1}$ belongs to $\OO_{K_{n+j}}[T]$, and satisfies $\val_K(R_n(0))  \geq M + \val_K(R_n'(0))$. By the theory of Newton polygons, $R_n(T)$ has a unique root of slope $\val_K(R_n(0)) - \val_K(R_n'(0)) \geq M$, and this root, which is $u_\tau^n - u_\tau^{n,j}$, therefore belongs to $K_{n+j}$. This implies that $u_\tau^n \in \OO_{K_{n+j}}$, which finishes the proof by induction on $n$.
\end{proof}

If $\tau$ satisfies the equivalent conditions of proposition \ref{algfil}, then we can write $u_\tau = f_\tau ( \phi_q^{-j_\tau}(u))$ for some $j_\tau \geq 0$ and $f_\tau(T) \in \OO_K \dcroc{T}$.

\begin{lemm}\label{ftnul}
We have $f_\tau(0)=0$, $f_\tau'(0) \neq 0$,  $P^\tau \circ f_\tau(T) = f_\tau \circ P(T)$ and $\Col_g^\tau \circ f_\tau(T) = f_\tau \circ \Col_g(T)$.
\end{lemm}

\begin{proof}
If $u_\tau = f_\tau(\phi_d^{-j}(u))$, then $P^\tau(u_\tau) = P^\tau \circ f_\tau (\phi_d^{-j}(u))$ and then $\phi_d(u_\tau) =  f_\tau \circ P (\phi_d^{-j}(u))$ so that $P^\tau \circ f_\tau(T) = f_\tau \circ P(T)$. Likewise, computing $g(u_\tau)$ in two ways shows that $\Col_g^\tau \circ f_\tau(T) = f_\tau \circ \Col_g(T)$. Evaluating  $P^\tau \circ f_\tau(T) = f_\tau \circ P(T)$ at $T=0$ gives $P^\tau(f_\tau(0)) = f_\tau(0)$ so that $f_\tau(0)$ is a root of $P^\tau(T) = T$. The theory of Newton polygons shows that those roots are $0$ and elements of valuation $0$. The latter case is excluded because $\theta \circ \phi_d^{-n}(u_\tau) = f_\tau ( u_{n+j} ) \in \MM_{K_\infty}$, so that $f_\tau(0) \in \MM_K$. We now prove that $f_\tau'(0) \neq 0$. Write $f(T) = f_k T^k + \bigO(T^{k+1})$ with $f_k \neq 0$. The fact that $P^\tau \circ f_\tau(T) = f_\tau \circ P(T)$ implies that $\tau(P'(0)) f_k = f_k P'(0)^k$ so that $\tau(P'(0)) = P'(0)^k$. Since $\val_K(P'(0)) > 0$, this implies that $k=1$.
\end{proof}

\begin{coro}\label{embssgr}
The set of those $\tau \in \Gal(K/\Qp)$ such that $r_\tau \geq 1$ forms a subgroup of $\Gal(K/\Qp)$, and if $F$ is the subfield of $K$ cut out by this subgroup, then $\eta(g) \in \OO_F^\times$. The weight $r_\tau$ is independent of $\tau \in \Gal(K/F)$.
\end{coro}

\begin{proof}
By proposition \ref{liftubar}, $\tau=\Id$ satisfies condition (iii) of proposition \ref{algfil} above, and therefore condition (i) as well, so that $r_{\Id} \geq 1$. If $\sigma, \tau$ satisfy condition (v) of ibid, then we can write $u_\sigma = f_\sigma(\phi_d^{-j_\sigma}(u))$ and $u_\tau = f_\tau(\phi_d^{-j_\tau}(u))$ so that $u_{\sigma \tau} = f_\tau^\sigma \circ f_\sigma (\phi_d^{-(j_\tau+j_\sigma)}(u))$ and therefore $\sigma \tau$ also satisfies condition (v). Since $\Gal(K/\Qp)$ is a finite group, these two facts imply that the set of $\tau \in \Gal(K/\Qp)$ such that $r_\tau \geq 1$ is a group.

By lemma \ref{ftnul}, we have $P^\tau \circ f_\tau(T) = f_\tau \circ P(T)$. This implies that $P'(0) \in \MM_F$ and also that $(P^\tau)^{\circ n} \circ f_\tau(T) = f_\tau \circ P^{\circ n}(T)$, so that 
\[ \frac{1}{P'(0)^n}(P^\tau)^{\circ n} \circ f_\tau(T) = \frac{1}{P'(0)^n} f_\tau \circ P^{\circ n}(T), \]
which implies by passing to the limit that $\Log^\tau \circ f_\tau(T) = f_\tau'(0) \cdot \Log(T)$. Since $\Col_g^\tau \circ f_\tau(T) = f_\tau \circ \Col_g(T)$, we have $g(\Log^\tau \circ f_\tau(u) ) = \tau(\eta(g)) \cdot (\Log^\tau \circ f_\tau(u))$. Moreover, $\Log^\tau \circ f_\tau(u) = f_\tau'(0) \cdot \Log(u)$, and therefore $\tau(\eta(g)) = \eta(g)$. This is true for every $\tau \in \Gal(K/F)$, so that $\eta(g) \in \OO_F^\times$. The fact that $\eta(g) \in \OO_F^\times$ implies that $r_\tau$ depends only on $\tau {\mid}_F$ and is therefore  independent of $\tau \in \Gal(K/F)$.
\end{proof}

\section{Local class field theory}
\label{lcft}

We now prove theorem D, and show how local class field theory allows us to derive theorem A from theorem D. We still assume that $K/\Qp$ is Galois for simplicity. Let $\lambda$ be a uniformizer of $\OO_K$ and let $K_\lambda$ denote the extension of $K$ attached to $\lambda$ by local class field theory. This extension is generated over $K$ by the torsion points of a Lubin-Tate formal group defined over $K$ and attached to $\lambda$ (see for instance \cite{LT} and \cite{LCFT}). Let $\chi^K_\lambda : \Gal(K_\lambda / K) \to \OO_K^\times$ denote the corresponding Lubin-Tate character. 

We still assume that the extension $K_\infty/K$ is Galois, so that it is an abelian totally ramified extension. This implies that there is a uniformizer $\lambda$ of $\OO_K$ such that $K_\infty \subset K_\lambda$. Let $\eta : \Gamma \to \OO_K^\times$ be the character constructed in \S\ref{chareta}.

\begin{prop}\label{etaexpl}
We have $\eta = \prod_{\tau \in \Gal(K/\Qp)} \tau(\chi^K_\lambda)^{r_\tau}$.
\end{prop}

\begin{proof}
The character $\eta : \Gamma \to \OO_K^\times$ is crystalline, and its weight at $\tau$ is $r_\tau$ by definition. The character $\eta_0 = \eta \cdot (\prod_{\tau \in \Gal(K/\Qp)}  \tau(\chi^K_\lambda)^{r_\tau} )^{-1}$ of $\Gal(K_\lambda/K)$ is therefore crystalline with weights $0$ at all embeddings, so that it is an unramified character of $\Gal(K_\lambda/K)$. Since $K_\lambda/K$ is totally ramified, we have $\eta_0=1$. 
\end{proof}

Proposition \ref{etaexpl} and corollary \ref{embssgr} imply the following, which is theorem D.

\begin{theo}\label{maineta}
There exists $F \subset K$ and $r \in \ZZ_{\geq 1}$ such that $\eta = \Nm_{K/F}(\chi^K_\lambda)^r$.
\end{theo}

We now show how this implies theorem A. If $u \in \OO_K^\times$, let $\mu^K_u$ denote the unramified character of $G_K$ that sends the Frobenius map of $k_K$ to $u$. If $F$ is a subfield of $K$, and $\Nm_{K/F}(\lambda) = \varpi^h u$ with $\varpi$ a uniformizer of $\OO_F$ and some $u \in \OO_F^\times$, then $\Nm_{K/F}(\chi^K_\lambda) = \chi^F_\varpi \cdot \mu^K_u$. 

\begin{prop}\label{chards}
Let $\dS$ be a relative Lubin-Tate group, attached to an extension $E/F$, and an element $\alpha = \varpi^h u \in \OO_F^\times$. The action of $\Gal(\Qpbar/E)$ on the torsion points of $\dS$ is given by $g(x) = [\chi^F_\varpi \cdot \mu^E_u(g)](x)$.
\end{prop}

\begin{proof}
See  \S 4 of \cite{Y08}.
\end{proof}

Let $F$ be the subfield of $K$ afforded by theorem \ref{maineta}, and let $E$ be the maximal unramified extension of $F$ contained in $K$.

\begin{theo}\label{dslt}
There exists a relative Lubin-Tate group $\dS$, relative to the extension $E/F$, such that if $K_\infty^{\dS}$ denotes the extension of $K$ generated by the torsion points of $\dS$, then $K_\infty \subset K_\infty^{\dS}$ and $K_\infty^{\dS} / K_\infty$ is a finite extension.
\end{theo}

\begin{proof}
Let $\lambda$ be a uniformizer of $K$ such that $K_\infty \subset K_\lambda$ and let $\pi = \Nm_{K/E}(\lambda)$ and $\alpha = \Nm_{K/F}(\lambda)$, so that $\pi$ is a uniformizer of $E$ and $\alpha = \Nm_{E/F}(\pi)$. Let $\dS$ be a relative Lubin-Tate group attached to $\alpha$, and let $K_\infty^{\dS}$ be the extension of $K$ generated by the torsion points of $\dS$. If $g \in \Gal(\Qpbar/K_\infty^{\dS})$, then $\Nm_{K/F}(\chi_K)(g)=1$ by proposition \ref{chards} and the observation preceding it, so that $\eta(g)=1$ by theorem \ref{maineta}. This implies that $K_\infty \subset K_\infty^{\dS}$. By Galois theory and theorem \ref{maineta},
\begin{enumerate}
\item $K_\infty^{\dS}$ is the field cut out by $\{ g \in G_K \mid \Nm_{K/F}(\chi_\lambda^K(g)) = 1\}$; 
\item $K_\infty$ is the field cut out by $\{ g \in G_K \mid \Nm_{K/F}(\chi_\lambda^K(g))^r = 1\}$.
\end{enumerate}
This implies that $K_\infty^{\dS} / K_\infty$ is a finite Galois extension, whose Galois group injects into $\{ x \in \OO_F^\times \mid x^r=1\}$.
\end{proof}

This proves theorem A. We conclude this {\S} with an example of a $\phi$-iterate extension that is Galois, corresponding to a polynomial $P(T) \in \Qp[T]$ such that $r=2$ and such that the extension $K_\infty^{\dS} / K_\infty$ is of degree $2$ in  the notation of theorems \ref{maineta} and  \ref{dslt}.

\begin{theo}\label{chebtheo}
Let $K=\QQ_3$, $P(T)=T^3+6T^2+9T$ and $u_0=-3$. The corresponding iterated extension $K_\infty$ is $\QQ_3(\mu_{3^\infty})^{\{\pm 1\} \subset \ZZ_3^\times}$, and $\eta=\chi_{\cyc}^2$.
\end{theo}

\begin{proof} 
For $k \geq 1$, let $C_k(T)$ denote the $k$-th Chebyshev polynomial, which is characterized by the fact that $C_k(\cos(\theta)) = \cos(k \theta)$. Let $P_k(T) = 2 ( C_k(T/2+1) - 1)$, so that $P_k(T)$ is a monic polynomial of degree $k$, and $P_k(2(\cos(\theta)-1)) = 2(\cos(k \theta)-1)$. Note that $P(T)=P_3(T)$ and that $u_0 = -3 = 2(\cos(2\pi/3)-1)$. The element $u_n$ is therefore a conjugate of $2(\cos(2\pi/3^{n+1})-1)$. This proves the fact that $K_\infty = \QQ_3(\mu_{3^\infty})^{\{\pm 1\} \subset \ZZ_3^\times}$

If $g \in G_{\QQ_3}$, then $g( 2(\cos(2\pi/3^n)-1) =  2(\cos(2\pi\chi_{\cyc}(g) /3^n)-1)$. This implies that $\Col_g(T) = P_k(T)$ if  $\chi_{\cyc}(g) = k \in \ZZ_{\geq 1}$. The formula for $\eta$ now follows from this, and the well-known fact that $C_k'(1)=k^2$ if $k \geq 1$. 
\end{proof}

We leave to the reader the generalization of this construction to other $p$ and other Lubin-Tate groups. The results of \S 2 of \cite{LMS} should be useful for this.

\providecommand{\bysame}{\leavevmode ---\ }
\providecommand{\og}{``}
\providecommand{\fg}{''}
\providecommand{\smfandname}{\&}
\providecommand{\smfedsname}{\'eds.}
\providecommand{\smfedname}{\'ed.}
\providecommand{\smfmastersthesisname}{M\'emoire}
\providecommand{\smfphdthesisname}{Th\`ese}

\end{document}